\documentclass[12pt]{amsart}
\usepackage{amssymb}

\begin{document}

\title[Holomorphic extension of smooth CR-mappings]
{Holomorphic extension of smooth CR-mappings between real-analytic
and real-algebraic CR-manifolds}
\author[F. Meylan, N. Mir, and D. Zaitsev]{Francine Meylan, Nordine Mir, and
Dmitri Zaitsev}
\address{F. Meylan : Institut de Math\'ematiques, Universit\'e de Fribourg, 1700
Perolles, Fribourg, Switzerland} \email{francine.meylan@unifr.ch}
\address{N. Mir : Universit\'e de Rouen, Laboratoire de Math\'ematiques
Rapha\"el Salem, UMR 6085 CNRS, 76821 Mont-Saint-Aignan Cedex,
France} \email{Nordine.Mir@univ-rouen.fr}
\address{D. Zaitsev : Mathematisches Institut, Eberhard-Karls-Universit\"at,  Auf der Morgenstelle 10
T\"ubingen, 72076 T\"ubingen, Germany}
\email{dmitri.zaitsev@uni-tuebingen.de}

\thanks{\noindent 2000 {{\em Mathematics Subject Classification.} 32V25, 32C16,
32H02, 32H40, 32V20. \\ The first author was partially
supported by Swiss NSF Grant 2100-063464.00/1}}

\def\Label#1{\label{#1}}
\def\1#1{\ov{#1}}
\def\2#1{\widetilde{#1}}
\def\6#1{\mathcal{#1}}
\def\4#1{\mathbb{#1}}
\def\3#1{\widehat{#1}}

\def\C{{\4C}}
\def\R{{\4R}}

\def\Re{{\sf Re}\,}
\def\Im{{\sf Im}\,}

\numberwithin{equation}{section}
\def\s{s}
\def\k{\kappa}
\def\ov{\overline}
\def\span{\text{\rm span}}
\def\ad{\text{\rm ad }}
\def\tr{\text{\rm tr}}
\def\xo {{x_0}}
\def\Rk{\text{\rm Rk\,}}
\def\sg{\sigma}
\def \emxy{E_{(M,M')}(X,Y)}
\def \semxy{\scrE_{(M,M')}(X,Y)}
\def \jkxy {J^k(X,Y)}
\def \gkxy {G^k(X,Y)}
\def \exy {E(X,Y)}
\def \sexy{\scrE(X,Y)}
\def \hn {holomorphically nondegenerate}
\def\hyp{hypersurface}
\def\prt#1{{\partial \over\partial #1}}
\def\det{{\text{\rm det}}}
\def\wob{{w\over B(z)}}
\def\co{\chi_1}
\def\po{p_0}
\def\fb {\bar f}
\def\gb {\bar g}
\def\Fb {\ov F}
\def\Gb {\ov G}
\def\Hb {\ov H}
\def\zb {\bar z}
\def\wb {\bar w}
\def \qb {\bar Q}
\def \t {\tau}
\def\z{\chi}
\def\w{\tau}
\def\Z{\zeta}
\def\phi{\varphi}
\def\eps{\varepsilon}

\def \T {\theta}
\def \Th {\Theta}
\def \L {\Lambda}
\def\b {\beta}
\def\a {\alpha}
\def\o {\omega}
\def\l {\lambda}

\def \im{\text{\rm Im }}
\def \re{\text{\rm Re }}
\def \Char{\text{\rm Char }}
\def \supp{\text{\rm supp }}
\def \codim{\text{\rm codim }}
\def \Ht{\text{\rm ht }}
\def \Dt{\text{\rm dt }}
\def \hO{\widehat{\mathcal O}}
\def \cl{\text{\rm cl }}
\def \bR{\mathbb R}
\def \bS{\mathbb S}
\def \bK{\mathbb K}
\def \bD{\mathbb D}
\def \bC{\mathbb C}
\def \C{\mathbb C}
\def \N{\mathbb N}
\def \bL{\mathbb L}
\def \bZ{\mathbb Z}
\def \bN{\mathbb N}
\def \scrF{\mathcal F}
\def \scrK{\mathcal K}
\def \mc #1 {\mathcal {#1}}
\def \scrM{\mathcal M}
\def \cR{\mathcal R}
\def \scrJ{\mathcal J}
\def \scrA{\mathcal A}
\def \scrO{\mathcal O}
\def \scrV{\mathcal V}
\def \scrL{\mathcal L}
\def \scrE{\mathcal E}
\def \hol{\text{\rm hol}}
\def \aut{\text{\rm aut}}
\def \Aut{\text{\rm Aut}}
\def \J{\text{\rm Jac}}
\def\jet#1#2{J^{#1}_{#2}}
\def\gp#1{G^{#1}}
\def\gpo{\gp {2k_0}_0}
\def\emmp {\scrF(M,p;M',p')}
\def\rk{\text{\rm rk\,}}
\def\Orb{\text{\rm Orb\,}}
\def\Exp{\text{\rm Exp\,}}
\def\Span{\text{\rm span\,}}
\def\d{\partial}
\def\D{\3J}
\def\pr{{\rm pr}}

\def\dbl{[\hskip -1pt [}
\def\dbr{]\hskip -1pt]}

\def \CZZ {\C \dbl Z,\zeta \dbr}
\def \D{\text{\rm Der}\,}
\def \Rk{\text{\rm Rk}\,}
\def \ima{\text{\rm im}\,}
\def \I {\mathcal I}

\newtheorem{Thm}{Theorem}[section]
\newtheorem{Def}[Thm]{Definition}
\newtheorem{Cor}[Thm]{Corollary}
\newtheorem{Pro}[Thm]{Proposition}
\newtheorem{Lem}[Thm]{Lemma}
\newtheorem{Rem}[Thm]{Remark}
\theoremstyle{definition}\newtheorem{Exa}[Thm]{Example}


\keywords{Reflection principle, CR-map, real-analytic
CR-submanifold, real-algebraic CR-submanifold, holomorphic
extension}

\maketitle

\section{Introduction and results}\Label{int}

The classical Schwarz reflection principle states that a
continuous map $f$ between real-analytic curves $M$ and $M'$ in
$\C$ that locally extends holomorphically to one side of $M$,
extends also holomorphically to a neighborhood of $M$ in $\C$. It
is well-known that the higher-dimensional analog of this statement
for maps $f\colon M\to M'$ between real-analytic CR-submanifolds
$M\subset\C^N$ and $M'\subset\C^{N'}$ does not hold without
additional assumptions (unless $M$ and $M'$ are totally real). In
this paper, we assume that $f$ is ${\mathcal C}^{\infty}$-smooth
and that the target $M'$ is {\em real-algebraic}, i.e.\ contained
in a real-algebraic subset of the same dimension. If $f$ is known
to be locally holomorphically extendible to one side of $M$ (when
$M$ is a hypersurface) or to a wedge with edge $M$ (when $M$ is a
generic submanifold of higher codimension), then $f$ automatically
satisfies the tangential Cauchy-Riemann equations, i.e.\ it is CR.
On the other hand, if $M$ is {\em minimal}, any CR-map $f\colon
M\to M'$ locally extends holomorphically to a wedge with edge $M$
by {\sc Tumanov}'s theorem \cite{Tu89} and hence, in that case,
the extension assumption can be replaced by assuming $f$ to be CR.

Local holomorphic extension of a CR-map $f\colon M\to M'$ may
clearly fail when $M'$ contains an (irreducible) complex-analytic
subvariety $E'$ of positive dimension and $f(M)\subset E'$.
Indeed, any nonextendible CR-function on $M$ composed  with a
nontrivial holomorphic map from a disc in $\C$ into $E'$ yields a
counterexample. Our first result shows that this is essentially
the only exception. Denote by $\6 E'$ the set  of all points
$p'\in M'$ through which there exist irreducible complex-analytic
subvarieties of $M'$ of positive dimension. We prove:

\begin{Thm}\Label{ncc0}
Let $M\subset \C^N$ and $M'\subset \C^{N'}$ be connected real-analytic and
real-algebraic CR-submanifolds respectively. Assume that $M$ is
minimal at a point $p\in M$. Then for any  ${\mathcal
C}^{\infty}$-smooth CR-map $f\colon M\to M'$, at least one of the
following conditions holds:
\begin{enumerate}
\item[(i)] $f$ extends holomorphically to a neighborhood of $p$ in $\C^N$;
\item[(ii)] $f$ sends a neighborhood of $p$ in $M$ into $\6E'$.
\end{enumerate}
\end{Thm}

If $M'$ is a real-analytic hypersurface, the set $\6E'$ consists
exactly of those points that are not of finite type in the sense
of {\sc D'Angelo} \cite{D1} (see {\sc Lempert} \cite{Le} for the
proof) and, in particular, $\6 E'$ is closed. The same fact also
holds if $M'$ is any real-analytic submanifold or even any
real-analytic subvariety (see \cite{D2}). However, in general,
$\6E'$ may not even be a real-analytic subset (see
Example~\ref{nonanalytic}). As an immediate consequence of Theorem
\ref{ncc0}, we have:

\begin{Cor}\Label{ncc}
Let $M\subset \C^N$ and $M'\subset \C^{N'}$ be connected real-analytic and
real-algebraic CR-submanifolds respectively. Assume that $M$ is
minimal at a point $p\in M$ and that all positive-dimensional
irreducible complex-analytic subvarieties in $M'$ are contained in
a fixed (complex-analytic) subvariety $V' \subset M'$. Then any
${\mathcal C}^{\infty}$-smooth CR-map $f\colon M\to M'$ that does
not send a neighborhood of $p$ in $M$ into $V'$ extends
holomorphically to a neighborhood of $p$ in $\C^N$.
\end{Cor}

 In view of an example due
to {\sc Ebenfelt} \cite{E96}, the minimality assumption on $M$ at
$p$ in Corollary \ref{ncc} cannot be replaced by the assumption
that $M$ is minimal somewhere.
On the other hand, if $M$ is also real-algebraic, this replacement
is possible:

\begin{Thm}\Label{nccalg}
Let $M\subset \C^N$ and $M'\subset \C^{N'}$ be connected
real-algebraic CR-submanifolds
 with $p\in M$ and let $V'\subset M'$ be as in Corollary~{\rm\ref{ncc}}.
Then the conclusion of Corollary~{\rm\ref{ncc}} holds provided $M$
is minimal somewhere.
\end{Thm}

In the setting of Theorem~\ref{nccalg}, any ${\mathcal
C}^{\infty}$-smooth CR-map $f\colon M\to M'$ that does not send a
neighborhood of $p$ in $M$ into $V'$ extends even algebraically to
a neighborhood of $p$ in $\C^N$ by a result of \cite{Z99} (see \S \ref{general}).
Since the subset $\6E'\subset M'$ is always closed, Corollary \ref{ncc} and Theorem \ref{nccalg}
imply:

\begin{Cor}\Label{stock}
Let $M\subset \C^N$ and $M'\subset \C^{N'}$ be connected
real-analytic and real-algebraic CR-submanifolds respectively.
Assume that $M$ is minimal at a point $p\in M$ and that $M'$ does
not contain any irreducible complex-analytic subvariety of positive dimension
through a point $p'\in M'$. Then any ${\mathcal
C}^{\infty}$-smooth CR-map $f\colon M\to M'$ with $f(p)=p'$
extends holomorphically to a neighborhood of $p$ in $\C^N$. The
same conclusion holds at a point $p\in M$ if $M$ is real-algebraic
and only somewhere minimal.
\end{Cor}

In the case when $M\subset \C^N$ is a real hypersurface,
 the first part of
Corollary~\ref{stock} is due to {\sc Pushnikov} \cite{Pu1,Pu2}
(see also  {\sc Coupet-Pinchuk-Sukhov} \cite{CPS00} for a similar
result). A prototype of a target real-algebraic CR-submanifold
with no nontrivial complex-analytic subvariety is given by the
unit sphere $\4S^{2N'-1}\subset \C^{N'}$. Even in that case,
Corollary \ref{stock} is new. Indeed, we have:

\begin{Cor}\Label{ncc-sphere}
Let $M\subset \C^N$ be a connected real-analytic CR-submanifold,
minimal at a point $p\in M$. Then any  ${\mathcal
C}^{\infty}$-smooth CR-map $f\colon M\to \4S^{2N'-1}$ extends
holomorphically to a neighborhood of $p$ in $\C^N$. The same
conclusion holds for any point $p\in M$ if $M$ is real-algebraic
and only somewhere minimal.
\end{Cor}

For $f$ of class $\6C^\infty$, Corollary~\ref{ncc-sphere} extends
results of {\sc Webster} \cite{W79}, {\sc Forstneri{\v c}} \cite{
Fo86,Fo89,Fo92}, {\sc Huang} \cite{H94} and {\sc
Baouendi-Huang-Rothschild} \cite{BHR}. (On the other hand, in
their setting, they prove holomorphic extension of $f$ of class
$\6C^k$ for appropriate $k$.)

If we restrict ourselves to {\em submersive} CR-maps (i.e.\  maps
for which the differential is surjective), a known obstruction to
their holomorphic extension is the {\em holomorphic degeneracy} of
the submanifolds. Recall that a real-analytic CR-submanifold $M$
is holomorphically degenerate (see {\sc Stanton} \cite{St}) at a
point $p\in M$ if there is a nontrivial holomorphic vector field
in a neighborhood of $p$ in $\C^N$ whose real and imaginary parts
are tangent to $M$. The existence of such a vector field and a
nonextendible CR-function on $M$ at $p$ yields nonextendible local
self CR-diffeomorphic maps of $M$ near $p$ (see \cite{BHR}). It is
known (see \cite{BER96}) that $M$ is holomorphically degenerate at
$p$ if and only if it is holomorphically degenerate everywhere on
the connected component of $p$. Our next result shows that for
source minimal CR-submanifolds, holomorphic degeneracy is
essentially the only obstruction for submersive CR-maps to be
holomorphically extendible.

\begin{Thm}\Label{hnd}
Let $M\subset \C^N$ and $M'\subset \C^{N'}$ be connected
real-analytic and real-algebraic CR-submanifolds respectively with
$p\in M$. Assume that $M$ is everywhere minimal and $M'$ is
holomorphically nondegenerate. Then any $\6C^\infty$-smooth CR-map
$f\colon M\to M'$ which is somewhere submersive
 extends holomorphically to a neighborhood of $p$ in $\C^N$.
\end{Thm}

In the case when $M\subset \C^N$ is a real hypersurface, a similar
result is due to \cite{CPS00}. Example \ref{crap} below shows that
the assumption that $M$ is {\em everywhere minimal} cannot be
replaced in Theorem \ref{hnd} by the weaker assumption that $M$
{\em is minimal at $p$}. On the other hand, if $M$ is
real-algebraic, a replacement with even weaker assumption on $M$
is possible:

\begin{Thm}\Label{alg-crap}
Let $M\subset \C^N$ and $M'\subset \C^{N'}$ be connected
real-algebraic CR-submanifolds with $p\in M$. Then the conclusion
of Theorem~{\rm\ref{hnd}} holds provided $M$ is somewhere minimal
and $M'$  is holomorphically nondegenerate.
\end{Thm}

In the setting of Theorem~\ref{alg-crap}, any ${\mathcal
C}^{\infty}$-smooth CR-map $f\colon M\to M'$ extends in fact
algebraically to a neighborhood of $p$ in $\C^N$ by a result of
\cite{Z99} (see \S \ref{general}). Theorem
\ref{alg-crap} extends a result of
{\sc Baouendi-Huang-Rothschild}
\cite{BHR} who obtained the same conclusion for
$M,M'\subset\C^N$ real-algebraic hypersurfaces and of {\sc
Kojcinovic} \cite{K} for $M,M'\subset\C^N$ generic submanifolds of
equal dimension. For further related results and history on the
analyticity problem for CR-mappings, the reader is referred to
\cite{Fo2,BER,Hua00}.

We shall derive the above results in \S \ref{general} from the
following statement that relates analyticity properties of a
smooth CR-map with geometric properties of its graph:

\begin{Thm}\Label{straight}
Let $M\subset \C^N$ and $M'\subset \C^{N'}$ be connected real-analytic and
real-algebraic CR-submanifolds respectively and $f\colon M\to M'$
a ${\mathcal C}^{\infty}$-smooth CR-map whose graph is denoted by
$\Gamma_f$. Assume that $M$ is minimal at a point $p\in M$ and
that $f$ does not extend holomorphically to any neighborhood of
$p$. Then there exists an integer $1\le n\le N'-1$ and a
real-analytic subset $A\subset M\times M'$ through $(p,f(p))$
containing a neighborhood $\Omega$ of $(p,f(p))$ in $\Gamma_f$ and
satisfying the following straightening property: for any point
$(q,f(q))$ in a dense open subset of $\Omega$, there exists a
neighborhood $U_q$ of $(q,f(q))$ in $\C^N\times\C^{N'}$ and a
holomorphic change of coordinates in $U_q$ of the form $(\tilde
z,\tilde z')=(z,\phi(z,z'))\in \C^N\times \C^{N'}$ such that
\begin{equation}\Label{real}
A\cap U_q = \{(z,z')\in U_q : z\in M, \; \tilde z'_{n+1}=\cdots=\tilde z'_{N'}=0\}.
\end{equation}
\end{Thm}

Theorem \ref{straight} will follow from the more general Theorem~\ref{straightplus},
where, similarly to \cite{CPS00},
the target $M'\subset\C^{N'}$ is assumed to be a real-algebraic subset and
an estimate for the number $n$ (in Theorem \ref{straight}) is given.
Our approach  follows partially the techniques initiated in
 \cite{Pu1,Pu2} and further extended in \cite{CPS00} in the case when $M$ is a hypersurface.
A crucial point in the proof of Theorem~\ref{straight} consists of showing
(after possible shrinking $M$ around $p$)
that near a generic point of the graph
$\Gamma_f$, the intersection of $M\times\C^{N'}$ with
the local Zariski closure of $\Gamma_f$ at $(p,f(p))$ (see \S \ref{Zar} for the definition)
is contained in $M\times M'$ (see Theorem
\ref{straightplus} and Proposition \ref{bourde}).
Here the arguments of \cite{Pu1,Pu2,CPS00} (e.g.\ the proof of Proposition 4.1 in \cite{CPS00})
do not seem to apply in our case
 and we have to proceed differently.
In \S \ref{mer} we give preliminary results based on a meromorphic extension
property obtained by the authors in the previous work \cite{MMZ1}.
\S\ref{Zar}--\ref{embed} are devoted to the proof of Theorem~\ref{straight}.

 {\bf Acknowledgment.} The authors would like to thank {\sc H.-M. Maire} for
his interest to the paper and many helpful discussions.

\section{Preliminaries and examples}

\subsection{CR-submanifolds and CR-maps}
A real submanifold $M\subset \C^N$ is called a CR-submanifold if
the dimension of the complex tangent space $T^c_pM:=T_pM\cap
iT_pM$ is independent of $p\in M$; $M$ is called {\em generic} if
for any point $p\in M$, one has $T_pM + iT_pM= T_p\C^N$. For a
CR-submanifold $M$ we write $T^{0,1}M:=T^{0,1}\C^N\cap \C TM$,
where $T^{0,1}\C^N$ is the bundle of $(0,1)$ tangent vectors in
$\C^N$. A function $h\colon M\to \C^{N'}$ of class ${\mathcal
C}^1$ is called a CR-function if for any section $L$ of the
CR-bundle, $Lf=0$. If $h$ is merely continuous, $h$ is still
called CR if it is annihilated by all vector fields $L$ as above
in the sense of distributions. A continuous map $f\colon M\to M'$
between CR-submanifolds $M\subset \C^N$ and $M'\subset \C^{N'}$ is
called a CR-map if all its components are CR-functions.

A CR-submanifold $M\subset \C^N$ is called {\em minimal} (in the
sense of {\sc Tumanov}) at a point $p\in M$ if there is no real
submanifold $S\subset M$ through $p$ with ${\rm dim}\, S<{\rm
dim}\, M$ and such that $T_q^cM\subset T_qS$, for all $q\in S$. It
is well-known that if $M$ is real-analytic, the minimality
condition of $M$ is equivalent to the finite type condition in the
sense of {\sc Kohn} and {\sc Bloom-Graham} (see \cite{BER}).

 A real (resp.\ complex) submanifold $M\subset \C^N$ is real-algebraic (resp.\ algebraic) if it is contained in a
 real-algebraic (resp.\ complex-algebraic) subvariety with the same real (resp.\ complex) dimension as
 that of $M$. A map $f\colon M\to M'$ between real submanifolds $M\subset \C^N$ and $M'\subset \C^{N'}$ is real-algebraic if its graph
 $\Gamma_f:=\{(z,f(z)):z\in M\}$ is a real-algebraic submanifold of $\C^N\times \C^{N'}$. Similarly, a holomorphic map between
 open subsets $\Omega\subset \C^N$ and $\Omega' \subset \C^{N'}$ is called algebraic if its graph is a complex-algebraic submanifold of
 $\Omega\times \Omega'$.

\subsection{Examples}

The following example shows that, even if $M'\subset \C^{N'}$ is a
real-analytic hypersurface, the subset $\6E'\subset M'$ of all
points that are not of finite D'Angelo type is not real-analytic
in general.

\begin{Exa}\Label{nonanalytic}
Consider the tube real-analytic hypersurface $M'\subset\C^4$ given by
\begin{equation}\Label{}
(\Re z_1)^2 - (\Re z_2)^2 + (\Re z_3)^2 = (\Re z_4)^3
\end{equation}
near the point $(1,1,0,0)\in M'$. We claim that the subset
$\6E'\subset M'$ is given by $\Re z_4\ge 0$ and is therefore not
analytic. Indeed, every intersection of $M'$ with $\{z_4={\rm
const},\, \Re z_4\ge 0\}$ contains complex lines through each
point and is hence everywhere of {\sc D'Angelo} infinite type. On
the other hand, if $\Re z_4<0$, the coordinate $\Re z_2$ can be
expressed as a strictly convex function of the other coordinates.
Therefore, $M'$ is strictly pseudoconvex at each such point and
thus of {\sc D'Angelo} finite type.
\end{Exa}

The following example shows that a somewhere  submersive $\6
C^{\infty}$-smooth CR-map $f\colon M\to M'$ between connected
real-analytic hypersurfaces in $\C^2$ can be real-analytic on some
connected component of the set of minimal points of $M$ and not
real-analytic in another component. In particular, the assumption
of Theorem~\ref{hnd} that $M$ is everywhere minimal cannot be
replaced by the weaker assumption that $M$ is minimal  at $p$.

\begin{Exa}\Label{crap}
As in {\sc Ebenfelt}'s example \cite{E96},
let $M,M'\subset\C^2$ be connected real-analytic hypersurfaces through $0$ given respectively by
$$\Im w= \theta(\arctan |z|^2, \Re w),  \quad \Im w= (\Re w) |z|^2, $$
where $t=\theta(\xi,s)$ is the unique solution of the algebraic equation $\xi(t^2+s^2) - t = 0$
 with $\theta(0,0)=0$ given by the implicit function theorem.
Note that $M$ and $M'$ are minimal
precisely outside the complex line $\{w=0\}$ and
that $M'$ is real-algebraic, but $M$ is not.
For every  $\6C^\infty$-smooth CR-function $\phi$ on $M$, define a map $f_\phi\colon M\to \C^2$ by
\begin{equation}\Label{}
f_\phi(z,w):=
\begin{cases}
(z,\,0) &       \Re w = 0 \\
(z,\, e^{-1/w}) & \Re w > 0 \\
(z+\phi(z,w)\;e^{1/w},\,0) & \Re w < 0.
\end{cases}
\end{equation}
By similar arguments as in \cite{E96} it follows that $f_\phi$ is
always a $\6C^\infty$-smooth CR-map sending $M$ into $M'$. Suppose
we can choose $\phi$ not holomorphically extendible to any
neighborhood in $\C^2$ of a fixed minimal point $p_0=(z_0,w_0)\in
M$ with $\Re w_0<0$. Then it is easy to see that $f_\phi$ is
somewhere submersive but does not extend holomorphically to any
neighborhood of the minimal point $p_0\in M$.

To show that the above choice of $p_0$ and $\phi$ is possible,
observe that $\theta$ can be factored as
$\theta(\xi,s)=s^2\xi (1+\2\theta(\xi,s))$ with $\2\theta$ analytic and vanishing at the origin.
Hence $\Im w\ge 0$ for every sufficiently small $(z,w)\in M$.
Then, for any real sufficiently small $x_0\ne 0$, the point $p_0:=(0,x_0)\in M$ is minimal
and a suitable branch of $e^{-1/(w-x_0)^{1/3}}$ extends to a
$\6C^\infty$-smooth CR-function $\phi$ on $M$
that is not holomorphically extendible to any neighborhood of $p_0$.
\end{Exa}

\section{A result on meromorphic extension and its applications}\Label{mer}

In what follows, for any subset $V\subset \C^k$, ${V}^*$ denotes
the set $\{\bar{z}:z\in V\}$ and, as usual, for any ring $A$, we
denote by $A[X]$, $X=(X_1,\ldots,X_s)$, the ring of polynomials in
$s$ indeterminates with coefficients in $A$. An important role in
the proof of Theorem \ref{straight} will be played by the
following meromorphic extension result from \cite[Theorem
2.6]{MMZ1}.

\begin{Thm}\label{NME}
Let $\Omega \subset\C^N$, $V\subset\C^k$ be open subsets,
$M\subset \Omega$ a connected generic real-analytic submanifold,
$G\colon M\to V$ a continuous CR-function and
$\Phi,\Psi\colon {V}^*\times \Omega \to\C$  holomorphic functions.
Assume that $M$ is minimal at every point
and that there exists a nonempty open subset of $M$ where $\Psi(\1{G(z)},z)$ does not vanish
and where the quotient
$$H(z):=\frac{\Phi (\1{G(z)},z)}{\Psi (\1{G(z)},z)}$$
is CR. Then $\Psi (\1{G(z)},z)$ does not vanish on a dense open subset $\widetilde{M}\subset M$
and $H$ extends from $\widetilde M$ meromorphically to a neighborhood of $M$ in $\C^N$.
\end{Thm}

Results in the spirit of Theorem \ref{NME}
have been important steps in proving regularity results for CR-mappings
(see e.g.\ \cite{Pu1, Pu2, BHR, CPS99, CPS00, MMZ1}).

For a generic real-analytic submanifold $M\subset \C^N$, denote by
${\mathcal C}^{\infty}(M)$ the ring of ${\mathcal
C}^{\infty}$-smooth functions on $M$, by ${\mathcal O}(M)$ the
ring of restrictions of holomorphic functions to $M$ and by
${\mathcal O}_p(M)$ the corresponding ring of germs at a point
$p\in M$. Similarly to \cite{CPS99} (see also
\cite{Pu1,Pu2,CPS00,MMZ1}), define a subring ${\6A}(M)\subset
{\mathcal C}^{\infty}(M)$ as follows: {\em a function $\eta \in
{\mathcal C}^{\infty}(M)$ belongs to ${\6A}(M)$ if and only if,
near every point $p\in M$, it can be written in the form $\eta (z)
\equiv \Phi (\overline{G(z)},z)$, where $G$ is a $\C^k$-valued
${\mathcal C}^{\infty}$-smooth CR-function in a neighborhood of
$p$ in $M$ for some $k$ and $\Phi$ is a holomorphic function in a
neighborhood of $(\1{G(p)},p)$ in $\C^k\times\C^N$}. Note that the
ring  ${\mathcal C}^{\omega}(M)$ of all real-analytic functions on
$M$ is a subring of ${\6A}(M)$. We have the following known
properties (see e.g.\ \cite{MMZ1}):

\begin{Lem}\Label{opht}
Let $M\subset \C^N$ be a connected generic real-analytic submanifold that is minimal at every point.
Then  for any $u\in {\6A}(M)$ the following hold:
\begin{enumerate}
\item[(i)] if $u$ vanishes on a nonempty open subset of $M$, then it vanishes identically on $M$;
\item[(ii)] if $L$ is a real-analytic $(0,1)$ vector field on $M$, then $Lu\in \6A(M)$.
\end{enumerate}
\end{Lem}

The following proposition is a consequence of Theorem \ref{NME} and will be essential for the proof of Theorem \ref{straight}.
In the proof we
follow the approach of \cite{Pu2} (see also \cite[Proposition 5.1]{CMS}).

\begin{Pro}\Label{manger}
Let $M\subset \C^N$ be a connected generic real-analytic submanifold
that is minimal at every point.
Let $F_1,\ldots,F_r$ be ${\mathcal C}^{\infty}$-smooth CR-functions on $M$
satisfying one of the following conditions:
\begin{enumerate}
\item [(i)] the restrictions of $F_1,\ldots,F_r$
to a nonempty open subset of $M$
satisfy a nontrivial polynomial identity
with coefficients in $\6A(M)$;
\item [(ii)] the restrictions of $F_1,\ldots,F_r, \1{F_1},\ldots,\1{F_r}$
to a nonempty open subset of $M$
satisfy a nontrivial polynomial
identity with coefficients in $\6C^\omega(M)$.
\end{enumerate}
Then for any point $q\in M$, the germs at $q$ of $F_1,\ldots,F_r$
satisfy a nontrivial polynomial identity
with coefficients in ${\6O}_q(M)$.
\end{Pro}

\begin{proof}
We first observe that, for the rest of the proof, we can assume
that the $(0,1)$ vector fields on $M$ are spanned by global
real-analytic vector fields on $M$. Indeed, suppose we have proved
Proposition~\ref{manger} under this additional assumption, then we
claim that Proposition~\ref{manger} follows from that case. For
this, for fixed $F_1,...,F_r$ as in Proposition~\ref{manger} (i)
(or (ii)), let $\Omega\subset M$ be
 the set of all points $q\in M$ for which the conclusion
 holds. Then $\Omega$ is clearly open. After shrinking $M$
 appropriately, we see that $\Omega \not = \emptyset$ by
 the above weaker supposed version of Proposition~\ref{manger}.
Analogously, shrinking $M$ around an accumulation point of $\Omega$,
we conclude that $\Omega$ is closed and therefore $\Omega=M$ as required.

Let now  $\6R(T)$ be a nontrivial polynomial
in $T=(T_1,\ldots,T_r)$ over $\6A(M)$ such that
\begin{equation}\Label{vanishing}
{\mathcal R}(F)|_{U} \equiv 0
\end{equation}
for some nonempty open subset $U\subset M$, where $F:=(F_1,\ldots,F_r)$.
We write $\6R(T)$ as a linear combination
\begin{equation}\Label{mon}
{\mathcal R}(T)=\sum_{j=1}^{l}\delta_jr_{j}(T),
\end{equation}
where each $\delta_j\ne 0$ is in $\6A(M)$ and $r_j$ is a monomial
in $T$. By Lemma \ref{opht}, each $\delta_j$ does not vanish on a
dense open subset of $M$. By shrinking $U$, we may assume that
$\delta_l$ does not vanish at every point of $U$. We prove the
desired conclusion by induction on the number $l$ of monomials
 in (\ref{mon}).
For $l=1$, (\ref{vanishing}) and (\ref{mon}) and the choice of $U$ imply that $r_1(F)|_U=0$.
Since $r_1$ is a monomial and each component of $\1{F}$ is in $\6A(M)$, it follows from Lemma~\ref{opht}
that $F_j=0$ for some $j$
which yields the required nontrivial polynomial identity
with coefficients in $\6O(M)$ (even in $\C$).

Suppose now that the desired conclusion holds for any polynomial $\6R$
whose  number of monomials is strictly less than $l$.
In view of (\ref{vanishing}) and (\ref{mon}) we have
\begin{equation}\Label{cms}
r_l(F)|_U +  \big(\sum_{j<l} \frac{\delta_j}{\delta_l} \, r_{j}(F)\big)|_U=0.
\end{equation}
Let $L$ be any global CR vector field on $M$ with real-analytic coefficients.
Applying $L$ to (\ref{cms}) and using the assumption that $F_j$ is CR for any $j$, we obtain
\begin{equation}\Label{obtain}
\Big(\sum_{j<l}L\big( \frac{\delta_j}{\delta_l}\big) r_{j}(F)\Big) |_U=0.
\end{equation}
By Lemma~\ref{opht} (ii), each coefficient $L(\delta_j/\delta_l)$
can be written as a ratio of two functions in ${\6A}(M)$. From
(\ref{obtain}), we are led to distinguish two cases. If  for some
$j\in \{1,\ldots,l-1\}$, $L(\delta_j/\delta_l)$ does not vanish
identically in $U$, then the required conclusion follows from the
induction hypothesis.

It remains to consider the case when
\begin{equation}\Label{cr-coeff}
L(\delta_j/\delta_l)=0,\quad {\rm in}\ U,
\end{equation}
for all $j$ and for all choices of $(0,1)$ vector field $L$.
Then (\ref{cr-coeff}) implies that each ratio $\delta_j/\delta_l$ is CR on $U$
by the assumption at the beginning of the proof.
Hence, by Theorem \ref{NME}, it follows that each $\delta_j/\delta_l$
extends meromorphically to a neighborhood of $M$ in $\C^N$
and therefore, (\ref{cms}) can be rewritten as
\begin{equation}\Label{cms1}
r_l(F)|_U +  \big(\sum_{j<l} m_j r_{j}(F)\big)|_U=0,
\end{equation}
with $m_1,\ldots,m_{l-1}$
being restrictions to $M$ of meromorphic functions in a neighborhood of $M$.
Since $M$ is connected and minimal everywhere, the identity
\begin{equation}\Label{id-claim}
r_l(F(z)) +  \sum_{j<l} m_j(z)r_{j}(F(z))=0
\end{equation}
holds for every $z\in M$ outside the set $S$ consisting of the poles of  the $m_j$'s.
This proves the desired conclusion under the assumption (i).

For the statement under the assumption (ii), consider a nontrivial polynomial
${\mathcal P}(T,\2{T})\in {\mathcal C}^{\omega}(M)[T,\2{T}]$
such that  ${\mathcal P}(F,\1{F})|_U=0$ for a non-empty open subset $U\subset M$.
We write
\begin{equation}\Label{flaw}
{\mathcal P}(T,\2{T})=\sum_{\nu \in \N^r, |\nu|\leq l}{\mathcal P}_{\nu}(\2{T}) T^{\nu},
\end{equation}
where each ${\mathcal P}_{\nu}(\2{T}) \in  {\mathcal
C}^{\omega}(M)[\2{T}]$ and at least one of the $\6P_\nu$'s is
nontrivial. If there exists $\nu_0\in \N^r$ such that ${\mathcal
P}_{\nu_0}(\1{F})$ is not zero in  the ring ${\6A}(M)$, then it
follows that the polynomial ${\mathcal Q}(T):={\mathcal
P}(T,\1{F})\in {\6A}(M)[T]$ is nontrivial and satisfies ${\mathcal
Q}(F)|_U=0$. Then the condition (i) is fulfilled and the required
conclusion is proved above.

It remains to consider the case
when ${\mathcal P}_{\nu}(\1{F})=0$ for any $\nu \in  \N^r$.
Fix any $\nu$ such that $\6P_\nu(\2T)$ is nontrivial. Let $\1 {\6 P_{\nu}}(T)$ denote the polynomial
in ${\mathcal C}^{\omega}(M)[T]$ obtained from $\6P_\nu$ by taking the complex conjugates of
its coefficients. Then $\1 {\6 P_{\nu}}(T)$ is a nontrivial polynomial in ${\mathcal A}(M)[T]$ and
satisfies $\1{\6 P_{\nu}}(F)=0$ on $M$. Here again, condition (i) is fulfilled and the desired conclusion
follows. The proof is complete.
\end{proof}

\section{Zariski closure of the graph of a CR-map}\Label{Zar}
Throughout this section, let  $M\subset \C^N$ be a real-analytic
generic submanifold, $p\in M$ a fixed point in $M$ and $f\colon
M\to \C^{N'}$ a ${\mathcal C}^{\infty}$-smooth CR-map. For $q\in
\C^N$, denote by ${\6O}_q$ the ring of germs at $q$ of holomorphic
functions in $\C^N$. The goal of this section is to define and
give some basic properties of the local Zariski closure of the
graph $\Gamma_f$ at $(p,f(p))$ over the ring ${\mathcal O}_p[z']$.

\subsection{Definition of the local Zariski closure}\Label{defzar}
For $M$, $f$ and $p$ as above, define the {\em $($local$)$ Zariski
closure} of $\Gamma_f$ at $(p,f(p))$ with respect to the ring
${\mathcal O}_p[z']$ as the germ $\6Z_f\subset\C^N\times\C^{N'}$
at $(p,f(p))$ of a complex-analytic set defined by the zero-set of
all elements in ${\mathcal O}_p[z']$ vanishing on $\Gamma_f$. Note
that since $\6Z_f$ contains the germ of the graph of $f$ through
$(p,f(p))$, it follows that ${\rm dim}_{\C}\, \6Z_f \geq N$. In
what follows, we shall denote by $\mu_p(f)$ the dimension of the
Zariski closure $\6Z_f$.

\begin{Rem}\Label{irred-remark}
{\rm Observe that if $M$ is furthermore assumed to be minimal at
$p$, all the components of the map $f$ extend to a wedge with edge
$M$ at $p$; in this case, it follows from unique continuation at
the edge that $\6Z_f$ is locally irreducible with respect to the
ring ${\mathcal O}_p[z']$.}
\end{Rem}
\subsection{Dimension of the local Zariski closure and transcendence degree}\Label{link}
In this section, we discuss a link between the dimension of the
Zariski closure $\mu_p(f)$ defined above and the notion of
transcendence degree considered in \cite{Pu1, Pu2, CMS,CPS00}. The
reader is referred to \cite{ZS} for basic notions from field
theory used here.

Since the ring $\6O_p(M)$ is an integral domain, one may consider its
quotient field that we denote by ${\scrM}_p(M)$.
Recall that, by a theorem of {\sc Tomassini} \cite{To66},
any germ in ${\mathcal O}_p(M)$ extends holomorphically to a neighborhood of $p$ in $\C^N$.
Hence an element belongs to ${\mathcal M}_p(M)$ if and only if it extends meromorphically  to a neighborhood of $p$ in $\C^N$.
 Note that if $M$ is {\em moreover assumed to be minimal at $p$}, it follows that the ring of germs at $p$ of
 ${\mathcal C}^{\infty}$-smooth CR-functions
on $M$ is an integral domain, which allows one to introduce its quotient field
containing ${\mathcal M}_p(M)$.
Therefore, for a generic submanifold
$M$ minimal at
$p$, one may consider the finitely generated field extension ${\mathcal
M}_p(M)(f_1,\ldots,f_{N'})$ over ${\mathcal M}_p(M)$ where
$f_1,\ldots,f_{N'}$ are the components of $f$ considered as germs at $p$.
(In the hypersurface case such a field extension has been studied  by {\sc Pushnikov} \cite{Pu1,Pu2}.)
The transcendence degree $m_p(f)$ of the above field extension is called the {\em
transcendence degree} of the CR-map $f$ at $p$ (see \cite{CMS,CPS00}).
We have the following standard relation between $m_p(f)$ and $\mu_p(f)$:

\begin{Lem}\Label{hp}
Let $M\subset \C^N$ be a generic real-analytic submanifold through some point $p\in M$ and
$f\colon M\to \C^{N'}$ a ${\mathcal C}^{\infty}$-smooth CR-map.
Assume that $M$ is minimal at $p$. Then $\mu_p(f)=N+m_p(f)$.
\end{Lem}

\begin{Rem}
{\rm The minimality of $M$ is needed to guarantee that ${\mathcal
M}_p(M)(f_1,\ldots,f_{N'})$ is a field so that the transcendence
degree is defined.}
\end{Rem}

The following well-known proposition shows the relevance of $\mu_p(f)$
 for the study of the holomorphic extension of $f$.

\begin{Pro}\Label{bof}
Let $M\subset \C^N$ be a generic real-analytic submanifold through a point $p$ and
$f\colon M\to \C^{N'}$ a ${\mathcal C}^{\infty}$-smooth CR-map.
Then the following are equivalent:
\begin{enumerate}
\item [(i)] $\mu_p(f)=N$;
\item [(ii)] $f$ is real-analytic near $p$.
\item [(iii)] $f$ extends holomorphically to a neighborhood of $p$ in $\C^N$.
\end{enumerate}
\end{Pro}

Proposition \ref{bof} is a consequence of theorems of {\sc Tomassini}
\cite{To66} and of {\sc Malgrange} \cite{Ma66}.

\section{Local geometry of the Zariski closure}\Label{NEW1}
\subsection{Preliminaries}\Label{momo}
We use the notation from \S \ref{Zar} and assume that $M$ is minimal at $p$ and that
\begin{equation}\Label{golf}
\mu_p(f)<N+N'
\end{equation}
holds. By shrinking $M$ around $p$ if necessary, we may assume
that $M$ is connected and minimal at all its points. In what
follows, for an open subset $\Omega \subset \C^k$, ${\mathcal
O}(\Omega)$ will denote the ring of holomorphic functions in
$\Omega$.

In \S \ref{Zar}, we saw that $\mu_p(f)\ge N$ and $m:=m_p(f)=\mu_p(f)-N$
coincides with the transcendence degree of
the field extension
${\mathcal M}_p(M)\subset {\mathcal M}_p(M)(f_1,\ldots,f_{N'})$, where $f=(f_1,\ldots,f_{N'})$.
This implies that there exist integers
$1\leq j_1<\ldots<j_m< N'$  such that
$f_{j_1},\ldots,f_{j_m}$ form a transcendence basis of
${\mathcal M}_p(M)(f)$ over
${\mathcal M}_p(M)$. After
renumbering the coordinates
$z':=(\zeta,w)\in \C^m\times \C^{N'-m}$ and setting
$m':=N'-m$, we
may assume that
\begin{equation}\Label{heure}
 f=(g,h)\in \C_{\zeta}^m\times \C_{w}^{m'},
\end{equation}
where $g=(g_1,\ldots,g_m)$ forms a transcendence basis of
${\mathcal M}_p(M)(f)$ over ${\mathcal M}_p(M)$.

Since the components of the germ at $p$ of the CR-map $h:M\to \C^{m'}$
are algebraically dependent over ${\mathcal M}_{p}(M)(g)$, there exist
 monic polynomials
$P_j(T)\in {\mathcal M}_p(M)(g)[T]$, $j=1,\ldots,m'$, such that if $h=(h_1,\ldots,h_{m'})$, then
\begin{equation}\Label{base}
P_j(h_j)=0,\ j=1,\ldots,m',\ {\rm in}\ {\mathcal M}_p(M)(f).
\end{equation}
As a consequence, there exist non-trivial polynomials $\widehat {P_j}(T)\in \6O_p(M)[g][T]$,
$j=1,\ldots,m'$, such that
\begin{equation}\Label{ase}
\widehat{P_j}(h_j)=0,\ j=1,\ldots,m'.
\end{equation}
For every $j=1,\ldots,m'$, we can write
\begin{equation}\Label{display1}
\widehat{P_j}(T)=\sum_{\nu \leq k_j}q_{j\nu}T^{\nu},
\end{equation}
where each $q_{j\nu}\in \6O_{p}(M)[g]$, $q_{jk_j}\not\equiv 0$ and $k_j\geq 1$. Since
each $q_{j\nu}$ is in $\6O_{p}(M)[g]$, we can also write
\begin{equation}\Label{display2}
q_{j\nu}=q_{j\nu}(z)=R_{j\nu}(z,g(z))
\end{equation}
where $R_{j\nu}(z,\zeta)\in \6O_p(M)[\zeta]$. Note that each $R_{j\nu}(z,\zeta)$ can also be viewed as
an element of ${\mathcal O}_p[\zeta]$.

Let $\Delta_p^N$ be a polydisc neighborhood of $p$ in $\C^N$ such that
the analogues of  (\ref{ase}) -- (\ref{display2}) hold with germs replaced
by their representatives in
$M\cap \Delta^N_p$ (denoted by the same letters).
Moreover, in view of Remark~\ref{irred-remark},
we may assume that the Zariski closure $\6Z_f$
can be represented by an irreducible (over the ring $\6O_p[z']$)
closed analytic subset of $\Delta_p^N\times\C^{N'}$
(also denoted by $\6Z_f$).
By shrinking $M$ we may also assume
that $M$ is contained in $\Delta_p^N$.
Hence we have
\begin{equation}\Label{inclusion}
\Gamma_f\subset \6Z_f\subset \Delta_p^N\times\C^{N'}.
\end{equation}
Define
\begin{equation}\Label{write}
\2{P_j}(z,\zeta;T):=\sum_{\nu=0}^{k_j}R_{j\nu}(z,\zeta)T^{\nu}\in {\mathcal O}(\Delta_{p}^N)[\zeta][T],\ j=1,\ldots,m'.
\end{equation}
It follows from (\ref{ase}) -- (\ref{display2}) that one has
\begin{equation}\Label{change}
\2{P_j}(z,g(z);h_j(z))=0,\ z\in M,\ j=1,\ldots,m'.
\end{equation}
Here each $R_{j\nu}(z,\zeta)\in {\mathcal O}(\Delta_p^N)[\zeta]$, $k_j\geq 1$, and
\begin{equation}\Label{1}
R_{jk_j}(z,g(z))\not \equiv 0,\quad  z\in M.
\end{equation}
Moreover, since ${\mathcal O}_p[\zeta][T]$ is a unique factorization domain (see e.g.\ \cite{ZS}) and
since $M$ is minimal at $p$, we may assume that the polynomials given by (\ref{write}) are irreducible.

Consider the complex-analytic variety $\6V_f\subset \C^N\times \C^{N'}$
through $(p,f(p))$ defined by
\begin{equation}\Label{subvariety}
\6V_f:=\{(z,\zeta,w)\in \Delta_p^N\times \C^{m}\times \C^{m'}:\2{P_j}(z,\zeta;w_j)=0,\ j=1,\ldots,m'\}.
\end{equation}
By (\ref{change}), $\6V_f$ contains the graph $\Gamma_f$
and hence the Zariski closure $\6Z_f$.
In fact, since by Lemma \ref{hp}, $\dim_{\C} \6Z_f=\mu_p(f)=N+m$,
it follows from the construction that $\6Z_f$ is the (unique)
irreducible component of $\6V_f$
(over $\6O_p[z']$)
containing $\Gamma_f$.
Note that $\6V_f$ is not irreducible in general and, moreover,
can have a dimension larger than $\mu_p(f)$.
(This may happen, e.g. if another component of $\6V_f$
is of higher dimension than $\mu_p(f)$).

For $j=1,\ldots,m'$, let $\2{D_j}(z,\zeta)\in {\mathcal O}(\Delta_p^N)[\zeta]$ be the
discriminant of the polynomial $\2{P_j}(z,\zeta;T)$ (with respect to $T$). Consider the complex-analytic set
\begin{equation}\Label{D}
\2{\mathcal D}:= \cup_{j=1}^{m'}\{(z,\zeta)\in \Delta_p^N \times \C^{m}: \2{D_j}(z,\zeta)=0\}.
\end{equation}
 By the irreducibility of each polynomial $\2{P_j}(z,\zeta;T)$, we have $\2{D_j}(z,\zeta)\not \equiv 0$ in $\Delta_p^N\times \C^m$, for $j=1,\ldots,m'$. Therefore from the algebraic independence of the components of the map $g$
over ${\mathcal M}_p(M)$, it follows that the graph of $g$ is not contained in $\2{\mathcal D}$, i.e.\ that for $z\in M$,
\begin{equation}\Label{discrim}
\2{D_j}(z,g(z))\not \equiv 0,\ {\rm for}\ j=1,\ldots,m'.
\end{equation}
By minimality of $M$ as before, the sets
$$\Sigma_j:=\{z\in M:\2{D_j}(z,g(z))=0\},\ j=1,\ldots,m',$$
are nowhere dense in $M$, and so is the set
\begin{equation}\Label{nude}
\Sigma:=\cup_{j=1}^{m'}\Sigma_j=\{z\in M: (z,g(z)) \in \2{\mathcal D}\}.
\end{equation}

\subsection{Description of $\6Z_f$ on a dense subset of the graph of $f$}\Label{curves}
By the implicit function theorem,
for any point $z_0\in M \setminus \Sigma$, there
exist polydisc neighborhoods of $z_0$, $g(z_0)$ and $h(z_0)$,
denoted by $\Delta_{z_0}^N
\subset \Delta_p^N\subset \C^N$, $\Delta^m_{g(z_0)}\subset \C^{m}$,
$\Delta_{h(z_0)}^{m'}\subset \C^{m'}$ respectively and a
holomorphic map
\begin{equation}\Label{vue}
\theta(z_0;\cdot):\Delta_{z_0}^N \times \Delta_{g(z_0)}^m\to \Delta_{h(z_0)}^{m'}
\end{equation} such that
for $(z,\zeta,w)\in  \Delta_{z_0}^N\times
\Delta_{g(z_0)}^m\times \Delta_{h(z_0)}^{m'}$,
\begin{equation}\Label{caract}
(z,\zeta,w)\in \6V_f \iff (z,\zeta,w)\in \6Z_f \iff  w=\theta(z_0;z,\zeta).
\end{equation}
Since $\Gamma_f\subset \6Z_f$ in view of (\ref{caract}), for every fixed
$z_0\in M\setminus \Sigma$, we have
\begin{equation}\Label{refer}
h(z)=\theta(z_0;z,g(z)),\ z\in M\cap \Delta_{z_0}^N.
\end{equation}
Let $Z_f\subset M\times \C^{N'}$ be the real-analytic subset given by
\begin{equation}\Label{coif}
Z_f:=\6Z_f\cap (M\times \C^{N'}),
\end{equation}
and, for every $z_0\in M\setminus \Sigma$, consider the real-analytic submanifold
$Z_f(z_0)\subset Z_f$ defined by setting
\begin{equation}\Label{coiffeur}
Z_f(z_0):=Z_f\cap (\Delta_{z_0}^N\times \Delta_{g(z_0)}^m\times \Delta_{h(z_0)}^{m'}).
\end{equation}
Note that $Z_f(z_0)$ contains the graph of $f$ over $M\cap \Delta_{z_0}^N$ and that, by making the holomorphic change of coordinates $(\2 z,\2 z')=(z,\phi (z,z'))\in \C^N\times \C^{N'}$ where
$\phi (z,z')=\phi (z,(\zeta,w)):=(\zeta,w-\theta (z_0;z,\zeta))$,
the submanifold $Z_f(z_0)$ is given in these new coordinates  by
\begin{equation}\Label{complex}
Z_f(z_0)=\{(\2 z,\2 z')\in \Delta_{z_0}^N\times \Delta_{g(z_0)}^m\times \C^{m'}:\2 z \in M,\ \2 z_{m+1}'=\ldots=\2 z_{N'}'=0\},
\end{equation}
where we write $\2 z'=(\2 z'_1,\ldots,\2 z'_{N'})$.

We summarize the above in the following proposition.

\begin{Pro}\Label{summa}
Let $M\subset \C^N$ be a generic real-analytic submanifold through a point $p\in M$ and $f:M\to \C^{N'}$
a ${\mathcal C}^{\infty}$-smooth CR-map. Let $\6Z_f$ be the local Zariski closure at $(p,f(p))$ of the graph of $f$
as defined in {\rm \S \ref{defzar}}. Assume that $M$ is minimal at $p$ and that $\mu_p(f)<N+N'$. Then after shrinking $M$ around
$p$, the following holds. For $z_0\in M\setminus \Sigma$, where $\Sigma$  is the nowhere dense open
subset of $M$ given by {\rm (\ref{nude})}, there exists a holomorphic change of coordinates near $(z_0,f(z_0))$ of the form
$(\2 z,\2 z')=(z,\phi (z,z'))\in \C^N\times \C^{N'}$ such that
the real-analytic subset $\6Z_f\cap (M\times \C^{N'})$ is given near $(z_0,f(z_0))$ by {\rm (\ref{complex})}, with $m=\mu_p(f)-N$.
\end{Pro}

For every $z_0\in M\setminus \Sigma$, denote by $\Omega_{z_0}$
the (unique) connected component of $(M\cap \Delta^N_{z_0})\times \Delta_{g(z_0)}^m$
passing through $(z_0,g(z_0))$.
Since $\Omega_{z_0}$
is connected, it makes sense to consider
the quotient field of the ring of real-analytic functions
${\mathcal C}^{\omega}(\Omega_{z_0})$ that we denote by $\6K(z_0)$.
Let
$$j\colon {\mathcal C}^{\omega}(M)[\zeta,\bar{\zeta}]\to {\mathcal C}^{\omega}(\Omega_{z_0})$$
 be the restriction map and
$$\6D:= {\rm Im}\,j \subset {\mathcal C}^{\omega}(\Omega_{z_0})$$
be the image of $j$.
Note that, since $\Omega_{z_0}$ is open in $M\times\C^m$, $j$ is an injective ring
homomorphism and hence, one can identify $\6D$ with
${\mathcal C}^{\omega}(M)[\zeta,\bar{\zeta}]$ via $j$.
Denote by $\6F$ the quotient field of $\6D$.
The field $\6F$ is naturally identified with the field
of all rational functions in $(\zeta,\bar\zeta)$
with coefficients that extend as real-analytic functions on $M$.
We have the field extension $\6F \subset \6K (z_0)$.

The following lemma, which will be needed for the proof of Theorem \ref{straight}, is  a direct
consequence of (\ref{caract}):

\begin{Lem}\Label{hope}
For every fixed $z_0\in M\setminus \Sigma$,
the restriction of the map $\theta (z_0;z,\zeta)$
$($given by {\rm (\ref{caract})}$)$ to $\Omega_{z_0}$
satisfies a nontrivial polynomial identity with coefficients in
$\6O(\Delta_p^N)[\zeta]$.
\end{Lem}

\section{Proof of Theorem \ref{straight}}\Label{embed}

With all the tools defined in \S \ref{Zar}--\S \ref{NEW1} at our
disposal,  we are now ready to prove the following statement from  which
Theorem \ref{straight} will follow.

\begin{Thm}\Label{straightplus}
Let $M\subset \C^N$ be a real-analytic generic submanifold through a point $p\in M$.
Let $f\colon M\to \C^{N'}$ be
a ${\mathcal C}^{\infty}$-smooth CR-map and $\6Z_f$ the local Zariski closure over
${\mathcal O}_p[z']$ at $(p,f(p))$ of $\Gamma_f$ as defined in {\rm \S \ref{defzar}}.
Suppose that $M$ is minimal at $p$ and $f$
maps $M$ into $M'$, where $M'$ is a proper real-algebraic subset of $\C^{N'}$.
Then, shrinking $M$ around $p$
and choosing an appropriate union $\2Z_f$
of local real-analytic irreducible components of $\6Z_f\cap (M\times \C^{N'})$ at $(p,f(p))$
if necessary, one has the following:
\begin{enumerate}
\item[(i)] $\mu_p(f)<N+N'$ for $\mu_p(f)=\dim\6Z_f$;
\item[(ii)] $\Gamma_f\subset\2Z_f\subset M\times M'$;
\item[(iii)] $\2Z_f$ satisfies the following straightening property:
for any point $q$ in a dense subset of $M$, there exists a neighborhood $U_q$ of $(q,f(q))$ in $\C^N\times\C^{N'}$
and a holomorphic change of coordinates in $U_q$ of the form $(\2z,\2z')=\Phi (z,z')=(z,\phi(z,z'))\in \C^N\times \C^{N'}$
such that
\begin{equation}\Label{bisreal}
\2Z_f\cap U_q = \{(z,z')\in U_q : z\in M, \; \2z'_{m+1}=\cdots=\2z'_{N'}=0\},
\end{equation}
where $m=\mu_p(f)-N$.
\end{enumerate}
\end{Thm}

For the proof we shall need the following result.

\begin{Pro}\Label{bourde}
Under the assumptions of Theorem~{\rm\ref{straightplus}},
shrinking $M$ around $p$ if necessary, one has the following:
\begin{enumerate}
\item[(i)] $\mu_p(f)<N+N'$;
\item[(ii)] For any point $z_0\in M\setminus \Sigma$, the real-analytic submanifold
$Z_f(z_0)$ is contained in $M\times M'$, where  $\Sigma$ is the nowhere dense subset of $M$
given by {\rm (\ref{nude})} and
$Z_f(z_0)\subset \6Z_f\cap (M\times \C^{N'})$ is  given by {\rm (\ref{coiffeur})}.
\end{enumerate}
\end{Pro}


\begin{proof}[Proof of Proposition {\rm \ref{bourde} (i)}]
We proceed by contradiction. Suppose that the dimension $\mu_p(f)$
 of the local Zariski closure is $N+N'$.
Since $M'$ is a proper real-algebraic subset of $\C^{N'}$,
there exists a nontrivial polynomial $\rho'(z',\1{z'})\in \C[z',\1{z'}]$
vanishing on $M'$.
Since $f$ maps $M$ into $M'$, we have
\begin{equation}\Label{goal}
\rho'(f(z),\overline{f(z)})=0
\end{equation}
for all $z\in M$.
It follows from Proposition \ref{manger} (ii)
(applied to $F:=f=(f_1,\ldots,f_{N'})$) that the germs at $p$ of
the components $f_1,\ldots,f_{N'}$ satisfy
a nontrivial polynomial identity with coefficients
in $\6O_p(M)$.
This contradicts the assumption $\mu_p(f)=N+N'$.
The proof of Proposition \ref{bourde} (i) is complete.
\end{proof}

In view of Proposition \ref{bourde} (i), we may now assume to be
in the setting of \S \ref{Zar}--\S \ref{NEW1}. Since $M'$ is
real-algebraic, it is given by
\begin{equation}\Label{orleans}
 M':=\{z'\in \C^{N'}:\rho'_1(z',\1{z'})=\ldots=\rho'_{l}(z',\1{z'})=0\},
 \end{equation}
 where each $\rho'_j(z',\1{z'})$, for $j=1,\ldots,l$, is a real-valued polynomial in $\C[z',\1{z'}]$.
For $j=1,\ldots,l$, $z_0\in M\setminus \Sigma$ and $(z,\zeta)\in \Omega_{z_0}$, define
\begin{equation}\Label{claim0}
r_j(z,\bar{z},\zeta,\bar{\zeta}):=\rho'_j(\zeta,\theta(z_0;z,\zeta),\bar{\zeta},\overline{
\theta (z_0;z,\zeta)})\in {\mathcal C}^{\omega}(\Omega_{z_0}),
\end{equation}
where $\theta (z_0;\cdot)\colon \Omega_{z_0} \to
\Delta_{h(z_0)}^{m'}$ is the restriction to $\Omega_{z_0}$ of the
holomorphic map given by (\ref{caract}) and $\Omega_{z_0}$ is the
open subset of $M\times \C^m$ given in \S \ref{curves}. We need
the following lemma.

\begin{Lem}\Label{claim}
For every $z_0\in M\setminus \Sigma$ and $j=1,\ldots,l$,
the real-analytic function $r_j$
satisfies a nontrivial polynomial identity on $\Omega_{z_0}$
with coefficients in ${\mathcal C}^{\omega}(M)[\zeta,\bar{\zeta}]$.
\end{Lem}

\begin{proof}
It follows from Lemma \ref{hope} that each component of the restriction to $\Omega_{z_0}$
of $\theta (z_0;\cdot)$, considered as an element of ${\mathcal C}^{\omega}(\Omega_{z_0})$,
is algebraic over the field $\6F$ defined in \S \ref{curves}.
Therefore, in view of the definition of $\6F$,
it is also the case for each component of the restriction to $\Omega_{z_0}$
of $\overline{\theta (z_0;\cdot)}$.
Since  for $j=1,\ldots,l$, each $\rho'_j(z',\1{z'})$ is a polynomial,
it follows from (\ref{claim0}) that each $r_j$ belongs
to the field generated by $\6F$ and the components of the restriction to $\Omega_{z_0}$
of the maps $\theta (z_0;\cdot)$ and $\overline{\theta (z_0;\cdot)}$.
Hence, by standard arguments  from field theory
(see e.g.\ \cite{ZS}),  each $r_j$ is also algebraic over $\6F$ for $j=1,\ldots,l$,
which gives the desired statement of the lemma.
\end{proof}

\begin{proof}[Proof of Proposition {\rm \ref{bourde} (ii)}.]
By contradiction, assume that there exists
$z_0\in M\setminus \Sigma$ such that
the real-analytic submanifold $Z_f(z_0)$ given by (\ref{coiffeur})
is not contained in $M\times M'$. In view of (\ref{caract}), (\ref{coif}), (\ref{coiffeur}) and (\ref{claim0}),
this means that
there exists $j_0\in  \{1,\ldots,l\}$ such that $r_{j_0}\not \equiv 0$ in
 $\Omega_{z_0}$. By Lemma
\ref{claim},
there exists a nontrivial polynomial $Q(z,\bar{z},
\zeta,\bar{\zeta};T)\in
{\mathcal C}^{\omega}(M)[\zeta,\bar{\zeta}][T]$ such that
\begin{equation}\Label{study}
Q(z,\bar{z},\zeta,\bar{\zeta};r_{j_0}(z,\bar{z},\zeta,\bar{\zeta}))\equiv 0,
\quad {\rm for}\ (z,\zeta)\in \Omega_{z_0}.
\end{equation}
Moreover, since $r_{j_0}$ does not vanish identically on $\Omega_{z_0}$
and $M$ is connected, we may choose $Q$ such that
\begin{equation}\Label{pou}
Q(z,\bar{z},\zeta,\bar{\zeta};0)\not \equiv 0\ {\rm for}\ (z,\zeta)\in M\times \C^m.
\end{equation}
Recall that we write $f=(g,h)$ as in (\ref{heure})
and that the graph of $g=(g_1,\ldots,g_m)$ over $M\cap \Delta_{z_0}^N$ is contained in $\Omega_{z_0}$.
Then (\ref{study}) implies that for $z\in M\cap \Delta_{z_0}^N$,
\begin{equation}\Label{hyp}
Q(z,\bar{z},g(z),\overline{g(z)};r_{j_0}(z,\bar{z},g(z),\overline{g(z)})) \equiv 0.
\end{equation}
But since $f$ maps $M$ into $M'$, we have for $j=1,\ldots,l$,
\begin{equation}\Label{maps}
\rho'_{j}(f(z),\overline{f(z)}) \equiv
\rho'_{j}(g(z),h(z),\overline{g(z)},\overline{h(z)}) \equiv 0, \quad z\in M.
\end{equation}
Therefore, in view of (\ref{refer}), (\ref{claim0}) and (\ref{maps}),
we obtain that for all $z\in M\cap \Delta_{z_0}^N$,
\begin{equation}\Label{push}
r_{j_0}(z,\bar{z},g(z),\overline{g(z)})  \equiv 0.
\end{equation}
From (\ref{hyp}) and (\ref{push}), we conclude that for all $z\in M\cap \Delta_{z_0}^N$,
\begin{equation}\Label{concl}
Q(z,\bar{z},g(z),\overline{g(z)};0) \equiv 0.
\end{equation}
In view of (\ref{pou}), condition (ii) in Proposition \ref{manger}
is satisfied for the components $g_1,\ldots,g_m$ of $g$
that are ${\mathcal C}^{\infty}$-smooth CR-functions on $M$.
By Proposition \ref{manger}, the germs at $p$ of
$g_1,\ldots,g_m$ satisfy a nontrivial polynomial
identity with coefficients in ${\6O}_p(M)$.
This contradicts the fact that $g_1,\ldots,g_m$  form
a transcendence basis of $\6M_p(M)(f)$ over $\6M_p(M)$ (see \S\ref{momo}).
The proof is complete.
\end{proof}

\begin{proof}[Proof of Theorem {\rm \ref{straightplus}}]
We shrink $M$ so that the conclusion of Proposition~\ref{bourde} holds.
Define $\2Z_f$ to be the union of those irreducible real-analytic components
of $\6Z_f\cap (M\times\C^{N'})$ that contain open pieces of $\Gamma_f$.
Then the conclusions (i) and (ii) of Theorem {\rm \ref{straightplus}} follow from Proposition~\ref{bourde}
and the straightening property (iii) follows from Proposition~\ref{summa}.
\end{proof}

\begin{proof}[Proof of Theorem {\rm \ref{straight}}] Without loss
of generality, we may assume that $M$ is generic. Since $f$ does
not extend holomorphically to any neighborhood of $p$ in $\C^N$,
we have $n:=\mu_p(f)-N>0$ by Proposition \ref{bof}. Then Theorem
\ref{straight} follows immediately from Theorem
\ref{straightplus}.
\end{proof}

\section{Proofs of Theorems \ref{ncc0}, \ref{nccalg}, \ref{hnd}, \ref{alg-crap}}\Label{general}

\begin{proof}[Proof of Theorem {\rm \ref{ncc0}}] We  need to prove that
if $f$ does not extend holomorphically to any neighborhood of $p$
in $\C^N$, then necessarily $f$ maps a neighborhood of $p$ in $M$
into $\6 E'$. By Theorem \ref{straight}, there exists a
neighborhood $U$ of $p$ in $M$ such that for all points $q$ in a
dense open subset of $U$, one has $f(q)\in \6 E'$. Since the set
$\6 E'$ is closed in $M'$ (see \S \ref{int}), it follows that
$f(U)\subset \6 E'$. This completes the proof of Theorem
\ref{ncc0}.
\end{proof}

\begin{proof}[Proof of Theorem {\rm \ref{nccalg}}] We may assume
that $M$ is generic. Since $M$ is real-algebraic, connected and
minimal somewhere, it is minimal outside a proper real-algebraic
subset $S$. In view of Corollary \ref{ncc}, we may assume that
$p\in S$.

If $W$ is a connected component of $M\setminus S$, then we claim
that either $f$ is real-algebraic on $W$ or $f(W) \subset V'$.
Indeed, if $f(W)\not \subset V'$, then $f$
 extends holomorphically to a neighborhood in $\C^N$
 of some point $q\in W$ by Corollary \ref{ncc}. Therefore, it
 is real-algebraic by a result of the third author \cite{Z99},
i.e.\ every component of $f$ satisfies a nontrivial polynomial
identity in a neighborhood of $q$ in $M$. Then by
 Tumanov's theorem and unique continuation, it follows that the
same polynomial identities for the components of $f$ hold
everywhere on $W$ and hence $f$ is real-algebraic on $W$.

By repeating the arguments from \cite[\S6]{DF78} one can show
that, near $p':=f(p)$, the set $V'$ (which may be empty) is
complex-algebraic, i.e.\ given by the vanishing of a vector-valued
holomorphic polynomial $P(z')$, $z'\in\C^{N'}$. Then, by the above
claim, $P\circ f$ is real-algebraic on each connected component of
$M\setminus S$. It is known  (see e.g.\ \cite{BR90}) that some
neighborhood of $p$ in $M$ intersects only finitely many connected
components of $M\setminus S$. Hence $P\circ f$ is real-algebraic
in a neighborhood of $p$ in $M$ and therefore, since $P\circ f$ is
$\6 C^{\infty}$-smooth, it is real-analytic by {\sc Malgrange}'s
theorem (see \cite{Ma66}).

If $f$ does not send a connected neighborhood of $p$ in $M$ into
$V'$, the real-analytic map $P\circ f$ does not vanish identically
on each of the components of $M\setminus S$ intersecting this
neighborhood. Hence, by the above claim, $f$ is real-algebraic on
each such component and therefore in a neighborhood of $p$. Then
the required holomorphic extension of $f$ at $p$ follows from {\sc
Malgrange}'s and {\sc Tomassini}'s theorems.
\end{proof}

For the proof of Theorem \ref{hnd}, it will be convenient to
derive the following corollary from Theorem \ref{straight}. For a
given smooth map $f\colon M\to M'$ between real submanifolds
$M\subset \C^N$ and $M'\subset \C^{N'}$ and a point $p\in M$, we
denote by ${\rm Rk}_p f$ the maximal rank of $f$ in a sufficiently
small neighborhood of $p$.

\begin{Cor}\Label{cor}
Under the assumptions of Theorem~{\rm\ref{straight}},
 there exists an integer $1\leq n\leq N'-1$
and a real-analytic submanifold $Y'\subset M'$ of dimension at
least ${\rm Rk}_pf$ arbitrarily close to $f(p)$ which is
biholomorphically equivalent to a product $Y\times \omega$, where
$Y\subset\C^{N'-n}$ is a real-analytic submanifold and $\omega$ is
an open subset in $\C^n$.
\end{Cor}

\begin{proof}[Proof of Corollary {\rm \ref{cor}}] Let $A,n$ be given by Theorem
\ref{straight}, $\pi'\colon \C^N_z\times \C^{N'}_{z'}\to
\C^{N'}_{z'}$ the natural projection and fix a sufficiently small
open neighborhood $U'$ of $f(p)$ in $\C^{N'}$. Choose  a point
$q\in M$ such that $f$ has rank ${\rm Rk}_pf$ at $q$ and $f(q)\in
U'$. Let $U_q,\phi$ be as in Theorem \ref{straight} and set
$\Phi(z,z'):=(z,\phi(z,z'))$. Shrinking $U_q$ if necessary, we may
assume that $f$ is of maximal rank at every point of $U_q$ and
that $\pi'(A\cap U_q)\subset U'$. Choose a point $(z_0,z'_0)\in
X_q:=A\cap U_q$ with $z'_0\in U'$, where $\pi'|_{X_q}:X_q\to
M'\cap U'$ has maximal rank $r'$. Since $\Gamma_f\cap U_q\subset
A$ by Theorem \ref{straight}, necessarily $r'\geq {\rm Rk}_pf$. By
the rank theorem, we may assume, after shrinking $X_q$ near
$(z_0,z_0')$ if necessary, that $Y':=\pi'(X_q)$ is an
$r'$-dimensional real-analytic submanifold in $U'$. By
Theorem~\ref{straight}, there exists a neighborhood $M_0\subset M$
of $z_0$ and an open subset $\omega\subset \C^n$ such that $\Phi
(X_q)=M_0\times \omega\times \{0\}$. Since $\pi'|_{X_q}$ and
$\pi'\circ \Phi^{-1}|_{\{z_0\}\times \omega \times \{0\}}$ are of
maximal rank, there exists a real-analytic submanifold $Y\subset
M_0$ such that, after shrinking $Y'$ and $\omega$ if necessary,
$\pi'\circ \Phi^{-1}$ defines a real-analytic diffeomorphism
between $Y\times \omega \times \{0\}$ and $Y'$. Since $\pi'$ and
$\Phi^{-1}$ are holomorphic in their respective ambient spaces, we
obtain the desired result.
\end{proof}

\begin{proof}[Proof of Theorem {\rm \ref{hnd}}]
Let $M$, $M'$ and $f$ be as in Theorem \ref{hnd}. Then $f$ is
submersive at a point $p_0\in M$. If $f$ were not holomorphically
extendible to a neighborhood of $p_0$ in $\C^N$, there would exist
an open holomorphically degenerate submanifold $Y'\subset M'$ by
Corollary \ref{cor}. This would contradict the assumption that
$M'$ is holomorphically nondegenerate. Hence $f$ is real-analytic
at $p_0$. Now define $\Omega\subset M$ to be the maximal connected
open subset containing $p_0$ where $f$ is real-analytic. Then $f$
is submersive on a dense subset of $\Omega$. We claim that
$\Omega=M$. Otherwise there would exist $p\in \1\Omega$ where $f$
is not real-analytic that would contradict Corollary \ref{cor} as
before. Hence $f$ is real-analytic everywhere on $M$ and the proof
is complete.
\end{proof}

\begin{proof}[Proof of Theorem {\rm \ref{alg-crap}}]
As in the proof of Theorem {\rm \ref{nccalg}}, we may assume that
$M$ is generic and we let $S\subset M$ be the real-algebraic
subset of all nonminimal points. By assumption, $f$ is submersive
at a point $p_0\in M$ which can be assumed minimal without loss of
generality. By Theorem~\ref{hnd}, $f$ extends holomorphically to a
neighborhood in $\C^N$ of the connected component $W_0$ of $p_0$
in $M\setminus S$. Since both $M$ and $M'$ are real-algebraic, $f$
is real-algebraic on $W_0$ by a result of \cite{Z99}. The same
argument shows that, for every connected component $W$ of
$M\setminus S$, either $f$ is real-algebraic or it is nowhere
submersive on $W$.

As in the proof of Theorem~\ref{hnd}, define $\Omega\subset M$ to
be the maximal open connected subset containing $p_0$ where $f$ is
real-analytic. Then $f$ is submersive on a dense subset of
$\Omega$. Assume by contradiction that $f$ is not real-analytic
everywhere on $M$ and hence that there exists a point
$p_1\in\1\Omega$ where $f$ is not real-analytic. Fix local
real-algebraic coordinates in $M$ and $M'$ near $p_1$ and $f(p_1)$
respectively and denote by $\Delta$ a minor of the Jacobian matrix
of $f$ of the maximal size that does not vanish identically in any
neighborhood of $p_1$. By the first part of the proof, we conclude
that $\Delta$ is real-algebraic and hence real-analytic in a
connected neighborhood $U$ of $p_1$ in $M$. In particular, $f$ is
submersive on a dense subset of $U$. Hence $f$ is real-algebraic
on every component of $U\setminus S$ by the first part of the
proof again, and hence, by {\sc Malgrange}'s theorem, it follows
that $f$ is real-analytic near $p_1$, which is a contradiction.
This shows that $\Omega=M$ and hence concludes the proof of the
theorem.
\end{proof}

\end{document}